\documentstyle{article}
\begin{document}
\font\germ=eufm10
\def\ssl{\hbox{\germ sl}}
\def\slh{\widehat{\ssl_2}}
\makeatletter
\def\aaa{@}
\centerline{\Large\bf Polyhedral Realizations of Crystal Bases for}
\centerline{\Large\bf  Quantized Kac-Moody Algebras}
\vskip12pt
\begin{tabular}{@{}cc}
\hspace{-1cm} Toshiki NAKASHIMA & Andrei ZELEVINSKY\\
&\\
\hspace{0cm}Department of Mathematical Science,& Department of Mathematics\\  
\hspace{0cm}Faculty of Engineering Science, &Northeastern University\\ 
\hspace{0cm}Osaka University, Osaka 560, Japan & Boston, MA 02115, USA \\ 
\hspace{0cm}toshiki\aaa sigmath.es.osaka-u.ac.jp & ANDREI\aaa neu.edu
\end{tabular}
\makeatother

   \renewcommand{\labelenumi}{(\roman{enumi})}
  \font\germ=eufm10
  
\def\al{\alpha}
\def\beneme{\begin{enumerate}}
\def\beq{\begin{equation}}
\def\beqn{\begin{eqnarray}}
\def\beqnn{\begin{eqnarray*}}
\def\bigsl{{\hbox{\fontD \char'54}}}
  \def\cd{\cdots}
  \def\del{\delta}
  \def\Del{\Delta}
  \def\ei{e_i}
  \def\eit{\tilde{e}_i}
\def\eneme{\end{enumerate}}
  \def\ep{\epsilon}
  \def\eeq{\end{equation}}
  \def\eeqn{\end{eqnarray}}
  \def\eeqnn{\end{eqnarray*}}
  \def\fit{\tilde{f}_i}
  \def\ge{\hbox{\germ g}}
  \def\gl{\hbox{\germ gl}}
  \def\hom{{\hbox{Hom}}}
  \def\ify{\infty}
  \def\io{\iota}
  \def\kp{k^{(+)}}
  \def\km{k^{(-)}}
  \def\llra{\relbar\joinrel\relbar\joinrel\relbar\joinrel\rightarrow}
  \def\lan{\langle}
  \def\lar{\longrightarrow}
  \def\lm{\lambda}
  \def\Lm{\Lambda}
  \def\mapright#1{\smash{\mathop{\longrightarrow}\limits^{#1}}}
  \def\nd{\noindent}
  \def\nn{\nonumber}
  \def\ot{\otimes}
  \def\op{\oplus}
  \def\ovl{\overline}
  \def\qq{\qquad}
  \def\q{\quad}
  \def\qed{\hfill\framebox[2mm]{}}
  \def\QQ{\hbox{\bf Q}}
  \def\qi{q_i}
  \def\qii{q_i^{-1}}
  \def\ran{\rangle}
  \def\ssl{\hbox{\germ sl}}
  \def\slh{\widehat{\ssl_2}}
  \def\ti{t_i}
  \def\tii{t_i^{-1}}
  \def\til{\tilde}
  \def\tt{{\hbox{\germ{t}}}}
  \def\ttt{\hbox{\germ t}}
  \def\uq{U_q(\ge)}
  \def\uqm{U^-_q(\ge)}
  \def\uqmq{{U^-_q(\ge)}_{\bf Q}}
  \def\uqq{U^{\bf Q}_q(\ge)}
  \def\util{\widetilde\uq}
  \def\vep{\varepsilon}
  \def\vp{\varphi}
  \def\vpi{\varphi^{-1}}
  \def\wtil{\widetilde}
  \def\what{\widehat}
  \def\ZZ{\hbox{\bf Z}}

\section{Introduction}
\setcounter{equation}{0}
\renewcommand{\theequation}{\thesection.\arabic{equation}}

Since pioneering works of G.Lusztig and M.Kashiwara on special bases
for quantum groups, a lot of work has been done 
on the combinatorial structure of these bases.
Although Lusztig's canonical bases and Kashiwara's global crystal bases
were shown by Lusztig to coincide whenever both are defined, 
their constructions are quite different and lead to different
combinatorial parametrizations.
In this paper we will only discuss the basis in the $q$-deformation
$U^-_q(\ge)$ of the universal enveloping algebra of the nilpotent 
part of a Kac-Moody Lie algebra $\ge$; understanding this basis 
is an essential first step towards understanding the bases in 
all the integrable highest weight modules of $\ge$.
Lusztig's approach (see \cite{Lu} and references there) works 
especially well when $\ge$
is a semisimple Lie algebra of simply-laced type.
In this case, every reduced expression for the maximal element 
of the Weyl group gives rise to a bijective parametrization of 
the canonical basis in $U^-_q(\ge)$ 
by the semigroup $\ZZ^N_{\geq 0}$ of all $N$-tuples of non-negative integers, 
where $N$ is the number of positive roots of $\ge$.
These parametrizations were studied in Lusztig's papers and also in
\cite{BFZ}.

Kashiwara's construction \cite{K0} of the global crystal basis in $U^-_q(\ge)$ 
is more elementary and works for an arbitrary Kac-Moody algebra.
The price for this is that the parametrizing sets for the basis
in Kashiwara's approach are more complicated than just $\ZZ^N_{\geq 0}$. 
In the literature, one can find several kinds of combinatorial 
expressions for crystal bases. 
As shown in \cite{KN},\cite{N}, crystal bases of finite-dimensional 
simple modules for classical Lie algebras  
can be parametrized by Young tableaux and their analogues.
For affine Lie algebras, crystal bases of integrable highest weight 
modules can be expressed as infinite sequences of 
perfect crystals (\cite{KMN1},\cite{KMN2}) or extended Young diagrams
(\cite{JMMO}). 
More generaly, in \cite{L1}, \cite{L2} for any symmetrizable 
Kac-Moody algebra, crystal bases are realized in terms of certain polygonal
paths. 
Although this presentation is elegant, it is not very convenient 
for actual computations with the basis.

In this paper we study bijective parametrizations of the crystal basis 
for $U^-_q(\ge)$ by integer sequences 
satisfying certain linear inequalities. 
In more geometric terms, the basis vectors should be parametrized
by lattice points in some polyhedral convex cone; 
this is what we mean by ``polyhedral realizations" in the title 
of the paper. 
If $\ge$ is semisimple, then, similarly to Lusztig's parametrizations,
a polyhedral realization is naturally associated with 
every reduced expression for the maximal element 
of the Weyl group.
This can be done using Kashiwara's theory of tensor products of crystals,
or, equivalently, using the ``string parametrizations" studied in 
\cite{BZ1}, \cite{BZ2}.
In fact, such a realization makes sense for arbitrary Kac-Moody
algebras, where reduced expressions for the maximal element 
of the Weyl group are replaced by certain infinite sequences of indices.

In this paper we deal with the following problem: describe explicitly
complete systems of linear inequalities that define all
polyhedral realizations of the crystal basis in $U^-_q(\ge)$. 
Such a description was recently found by P. Littelmann (private communication)
for all semisimple Lie algebras, using a case-by-case analysis.
Littelmann only treats some specific choice of a reduced expression
(for the type $A_n$, the corresponding result was already obtained in
\cite{BZ1}). 
We would like to find a unified description of all polyhedral realizations
for an arbitrary Kac-Moody algebra.
For Kac-Moody algebras of rank 2, such a description was found 
by M.Kashiwara in \cite{K3}, Proposition 2.2.3 (this description is sharpened
in Theorem \ref{rank 2-thm} below). 
The main result of the present paper (Theorem \ref{main} below) is a 
generalization of this description to Kac-Moody algebras 
of an arbitrary rank. 
The answer is given in terms of certain piecewise-linear transformations
of $\ZZ^{\ify}$. 
Unfortunately, our main result is only proved under certain 
technical assumptions.
These assumptions are checked to be valid in many cases including the 
ordinary and affine Lie algebras of type $A$.
It is even conceivable that they are always satisfied
(see the discussion in Section 3 below).  
 
The paper is organized as follows. 
In section 2, we review crystals and their basic properties.
We also introduce our main object of study: 
the crystal $B(\ify)$ corresponding to $U^-_q(\ge)$. 
We then describe the Kashiwara embedding of $B(\ify)$
into the lattice $\ZZ^{\ify}$. 
Our main theorem is formulated and proved in Section 3:
this is a description of the image of the Kashiwara embedding.
In the rest of the paper (sections 4, 5, and 6), 
we apply the theorem to the cases
when $\ge$ is of rank 2, of type $A_n$, and of type $A^{(1)}_{n-1}$,
respectively.

In preparing this paper, we received the preprint 
``Crystal bases and Young Tableaux" by G.Cliff. 
Thit paper describes the image of the Kashiwara embedding 
for the types $A,B,C,D$ and some special reduced expressions.
The method used in the preprint is 
different from ours. It seems that the method cannot be applied 
to affine algebras or more general Kac-Moody algebras,
or even to other reduced expressions for classical Lie algebras.
In a forthcoming paper, polyhedral realizations will be described 
not only for $B(\ify)$ but also for the irreducible $\ge$-modules.

This work was partly done  during the stay of T.N. at 
Northeastern University. He is grateful to the colleagues 
for their kind hospitality.
The work of A.Z. was partially supported by the NSF grant 
DMS-9625511.

\section{Preliminaries on Crystals}
\setcounter{equation}{0}
\renewcommand{\theequation}{\thesection.\arabic{equation}}
\subsection{Definition of $U_q(\ge)$}

Let $\ge$ be 
a  symmetrizable Kac-Moody algebra over {\bf Q} 
with a Cartan subalgebra 
$\ttt$, a weight lattice $P \subset \ttt^*$, the set of simple roots 
$\{\al_i: i\in I\} \subset \ttt^*$, 
and the set of coroots $\{h_i: i\in I\} \subset \ttt$,   
where $I$ is a finite index set 
(see \cite{Kac} for the background on Kac-Moody algebras). 
Let $\lan h,\lm\ran$ be the pairing between $\ttt$ and $\ttt^*$, 
and $(\al, \beta)$ be an inner product on 
$\ttt^*$ such that $(\al_i,\al_i)\in 2{\bf Z}_{\geq 0}$ and 
$\lan h_i,\lm\ran=2(\al_i,\lm)/(\al_i,\al_i)$ 
for $\lm\in\ttt^*$.
Let $P^*=\{h\in \ttt: \lan h,P\ran\subset\ZZ\}$. 

As in \cite{K1}, we define the quantized enveloping algebra $\uq$ 
to be an associative 
$\QQ(q)$-algebra generated by the $e_i$, $f_i \,\, (i\in I)$, 
and $q^h \,\, (h\in P^*)$ 
satisfying the following relations:
\begin{eqnarray}
&&q^0=1, \q{\hbox{\rm and }}\q q^hq^{h'}=q^{h+h'},\\
&&q^he_iq^{-h}=q^{\lan h,\al_i\ran}e_i,\qq
q^hf_iq^{-h}=q^{-\lan h,\al_i\ran}f_i,\\
&&e_i f_j - f_j e_i = \del_{i,j}(t_i-t^{-1}_i)/(q_i-q^{-1}_i),\\
&&\sum_{k=1}^{1-{\lan h_i,\al_j\ran}}
(-1)^kx_i^{(k)}x_jx_i^{(1-{\lan h_i,\al_j\ran}-k)}=0,\qq(i\ne j)
\end{eqnarray}
where the symbol $x_i$ in (2.4) stands for $e_i$ or $f_i$, and we set 
$q_i=q^{(\al_i,\al_i)/2}$, 
$t_i=q_i^{h_i}$, $[l]_i=(q^l_i-q^{-l}_i)/(q_i-q_i^{-1})$, 
$[k]_i!=\prod_{l=1}^k [l]_i$,
and $x_i^{(k)}=x_i^k/[k]_i!$. 

It is well-known \cite{CP} that $\uq$ has a Hopf algebra structure 
with the comultiplication $\Del$ given by 
$$
\Del(q^h)=q^h\ot q^h,\q \Del(\ei)=\ei\ot\tii+1\ot\ei,\q
\Del(f_i)=f_i\ot 1+\ti\ot f_i,
$$
for any $i\in I$ and $h\in P^*$. 
By this comultiplication, the tensor product of $\uq$-modules has 
a $\uq$-module structure.

\subsection{Definition of crystals}

The following definition of a crystal is due to  M.Kashiwara \cite{K3,K4};
it is motivated by abstracting some combinatorial properties 
of crystal bases. 
In what follows we fix a finite index set $I$ and a 
weight lattice $P$ as above.

\newtheorem{df2}{Definition}[section]
\begin{df2}
\label{def2.1}
A {\it crystal} $B$ is a set endowed with the following maps:
\begin{eqnarray}
&& wt:B\lar P,\\
&&\vep_i:B\lar\ZZ\sqcup\{-\infty\},\q
  \vp_i:B\lar\ZZ\sqcup\{-\infty\} \q{\hbox{for}}\q i\in I,\\
&&\eit:B\lar B\sqcup\{0\}, 
\q\fit:B\lar B\sqcup\{0\}\q{\hbox{for}}\q i\in I.
\end{eqnarray}
Here $0$ is an 
ideal element which is not included in $B$.
These maps must satisfy the following axioms: 
for all $b$,$b_1$,$b_2\in B$, we have
\begin{eqnarray}
&&\vp_i(b)=\vep_i(b)+\lan h_i,wt(b)\ran,\\
&&wt(\eit b)=wt(b)+\al_i{\hbox{ if }}\eit b\in B,\\
&&wt(\fit b)=wt(b)-\al_i{\hbox{ if }}\fit b\in B,\\
&&\eit b_2=b_1 {\hbox{ if and only if }} \fit b_1=b_2,
\label{eeff}\\
&&{\hbox{if }}\vep_i(b)=-\infty,
  {\hbox{ then }}\eit b=\fit b=0.
\end{eqnarray}
\end{df2}

The above axioms allow us to make a crystal $B$ into 
a colored oriented graph with the set of colors $I$.
This means that each edge of the graph is labeled with some $i \in I$;
we write $b_1\mapright{i} b_2$ for an oriented edge from $b_1$ to $b_2$ 
labeled with $i$.

\begin{df2}
\label{c-gra}
The crystal graph of a crystal $B$ is 
a colored oriented graph given by
the rule : $b_1\mapright{i} b_2$ if and only if $b_2=\fit b_1$  
$(b_1,b_2\in B)$.
\end{df2}

\begin{df2}
\label{df:mor}
\begin{enumerate}
\item
Let $B_1$ and $B_2$ be crystals. 
A morphism of crystals $\psi:B_1\lar B_2$ 
is a map $\psi:B_1 \lar B_2\sqcup\{0\}$ 
satisfying the following conditions:
\begin{eqnarray}
&&\hspace{-30pt}wt(\psi(b)) = wt(b),\q \vep_i(\psi(b)) = \vep_i(b),\q 
\vp_i(\psi(b)) = \vp_i(b)
\label{wt}\\
&&{\hbox{if }}b\in B_1{\hbox{ and }}\psi(b)\in B_2,\nonumber\\
&&\hspace{-30pt}\psi(\eit b)
=\eit\psi(b){\hbox{ if }}b\in B_1{\hbox{ satisfies }}
 \psi(b)\neq0{\hbox{ and }}\psi(\eit b)\neq0,\\
&&\hspace{-30pt}\psi(\fit b)
=\fit\psi(b){\hbox{ if }}b\in B_1{\hbox{ satisfies }}
 \psi(b)\neq0{\hbox{ and }}\psi(\fit b)\neq0.
\end{eqnarray}
\item
A morphism of crystals $\psi:B_1\lar B_2$ 
is called {\it strict} if the 
map $\psi: B_1 \lar B_2\sqcup\{0\}$ 
commutes with all $\eit$ and $\fit$ 
An injective strict morphism is called an embedding of crystals. 
\end{enumerate} 
\end{df2}

For crystals $B_1$ and $B_2$, we define their tensor product 
$B_1\ot B_2$ as follows:
\begin{eqnarray}
&&B_1\ot B_2=\{b_1\ot b_2: b_1\in B_1 ,\, b_2\in B_2\},\\
&&wt(b_1\ot b_2)=wt(b_1)+wt(b_2),\\
&&\vep_i(b_1\ot b_2)={\hbox{max}}(\vep_i(b_1),
  \vep_i(b_2)-\lan h_i,wt(b_1)\ran),
\label{tensor-vep}\\
&&\vp_i(b_1\ot b_2)={\hbox{max}}(\vp_i(b_2),
  \vp_i(b_1)+\lan h_i,wt(b_2)\ran),
\label{tensor-vp}\\
&&\eit(b_1\ot b_2)=
\left\{
\begin{array}{ll}
\eit b_1\ot b_2 & {\mbox{ if }}\vp_i(b_1)\geq \vep_i(b_2)\\
b_1\ot\eit b_2  & {\mbox{ if }}\vp_i(b_1)< \vep_i(b_2),
\end{array}
\right.
\label{tensor-e}
\\
&&\fit(b_1\ot b_2)=
\left\{
\begin{array}{ll}
\fit b_1\ot b_2 & {\mbox{ if }}\vp_i(b_1)>\vep_i(b_2)\\
b_1\ot\fit b_2  & {\mbox{ if }}\vp_i(b_1)\leq \vep_i(b_2).
\label{tensor-f}
\end{array}
\right.
\end{eqnarray}
Here $b_1\ot b_2$ is just another notation for an ordered pair
$(b_1,b_2)$, and we set $b_1 \ot 0 = 0 \ot b_2=0$. 
Let ${\cal C}(I,P)$ be the category 
of crystals with the index set $I$ and the weight lattice $P$. 
Then $\ot$ is a functor 
from ${\cal C}(I,P)\times{\cal C}(I,P)$ 
to ${\cal C}(I,P)$ that makes ${\cal C}(I,P)$ a tensor category 
\cite{Kas}.
In particular, the tensor product of crystals is
associative: the crystals 
$(B_1\ot B_2)\ot B_3$ and $B_1\ot(B_2\ot B_3)$ are isomorphic via 
$(b_1\ot b_2)\ot b_3\leftrightarrow b_1\ot (b_2\ot b_3)$.

We conclude this subsection with an example of a crystal
that will be needed later.

\newtheorem{ex}[df2]{Example}
\begin{ex}
\label{Example:crystal}
For $i\in I$, the crystal $B_i:=\{(x)_i\,: \, x \in\ZZ\}$ is defined by
\beqnn
&& wt((x)_i)=x \al_i,\qq \vep_i((x)_i)=-x,\qq \vp_i((x)_i)=x,\\
&& \vep_j((x)_i)=-\infty,\qq \vp_j((x)_j)
   =-\infty \q {\rm for }\q j\ne i,\\
&& \eit (x)_i=(x+1)_i,\qq \fit(x)_i=(x-1)_i,\\
&& \til e_j(x)_i=\til f_j(x)_i=0\q {\rm for }\q j\ne i.
\eeqnn
\end{ex}

\subsection{Crystal base of $U^-_q(\ge)$ and the crystal $B(\ify)$}

In this subsection we introduce the crystal $B(\ify)$, our main 
object of study.
All the results below are due to M.Kashiwara \cite{K1}. 
Let $U^-_q(\ge)$ be the subalgebra of $\uq$ generated by 
$\{f_i\}_{i\in I}$.
By Lemma 3.4.1 in \cite{K1}, for any $u \in\uqm$ and $i\in I$, 
there exist unique 
$u', u'' \in \uqm$ such that 
\beq
e_i u - u e_i ={{t_i u'' - t^{-1}_i u'}\over{q_i-q^{-1}_i}}.
\label{eiP}
\eeq
We define the endomorphisms $e_i'$ and  $e''_i$ of $\uqm$
by setting $e_i'(u)=u'$ and $e_i''(u)=u''$. 
For any $i\in I$, we have the 
direct sum decomposition 
\beq
\uqm=\bigoplus_{k \geq 0} f^{(k)}_i{\rm Ker}\,e'_i\q.
\label{dec-uqm}
\eeq 
Using this, we can define the endomorphisms 
$\eit$ and $\fit$ of $\uqm$ 
by
\beq
\eit(f^{(k)}_iu)=f^{(k-1)}_iu,\,\, {\rm and}\,\,
\fit(f^{(k)}_iu)=f^{(k+1)}_iu\q {\rm for}\,\, u\in\,{\rm Ker }\,\,e_i'.
\label{eit-fit}
\eeq

Let $A\subset \QQ(q)$ be the subring of rational functions 
that are regular at $q=0$.
Let $L(\ify)$ be the left $A$-submodule of $\uqm$ generated 
by all the elements $\til f_{i_l}\cd \til f_{i_1}\cdot 1$
with $l \geq 0$ and $i_j \in I$.
Then $L(\ify)/qL(\ify)$ is a $\QQ$-vector space.
We define a subset $B(\ify) \subset L(\ify)/qL(\ify)$
to be the set of all non-zero elements of the form
$\til f_{i_l}\cdot\cdot \til f_{i_1}\cdot 1\,
{\rm mod}\,qL(\ify)$. 
The pair $(L(\ify),B(\ify))$ is called the  
{\it crystal base } of $\uqm$. 
It satisfies the following properties:
\begin{enumerate}
\item
$L(\ify)$ 
is a free $A$-submodule of $\uqm$, 
and $\uqm\cong \QQ(q)\ot_A L(\ify)$. 
\item
$B(\ify)$ 
is a basis of the $\QQ$-vector space $L(\ify)/qL(\ify)$.
\item
The endomorphisms $\eit$ and $\fit$ preserve $L(\ify)$, and so act on 
$L(\ify)/qL(\ify)$.
\item
For any $i \in I$, we have 
$\eit B(\ify)\subset B(\ify)\sqcup\{0\}$ 
and
$\fit B(\ify)\subset B(\ify)$. 
\item
For $u,v\in B(\ify)$, we have
$\fit u=v$ if and only if $\eit v=u$.
\end{enumerate}

We denote by $u_{\ify} \in B(\ify)$ the image of $1$ under the projection
$L(\ify)  \lar L(\ify)/qL(\ify)$, and define the weight function 
$wt:B(\ify)\lar P$ by 
$wt(b):=-(\al_{i_1}+\cd+\al_{i_l})$ for 
$b=\til f_{i_l}\cd\til f_{i_1}u_{\ify}$.
We define integer-valued functions $\vep_i$ and $\vp_i$ on 
$B(\ify)$ by 
$$\vep_i(b):={\rm max} \q \{k: \eit^k b\ne 0\}, \,\, 
\vp_i(b):=\lan h_i,wt(b)\ran +\vep_i(b).$$
An easy check shows that $B(\ify)$ equipped with the operators 
$\eit$ and $\fit$, and with the functions $wt$, $\vep_i$ and $\vp_i$ 
is a crystal.

\subsection{Kashiwara embeddings of $B(\ify)$}
Consider the additive group
\beq 
\ZZ^{\ify}
:=\{(\cd,x_k,\cd,x_2,x_1): x_k\in\ZZ
\,\,{\rm and}\,\,x_k=0\,\,{\rm for}\,\,k\gg 0\};
\label{uni-cone}
\eeq
we will denote by $\ZZ^{\ify}_{\geq 0} \subset \ZZ^{\ify}$ 
the subsemigroup of nonnegative sequences. 
To the rest of this section, we fix an infinite sequence of indices 
$\io=(\cd,i_k,\cd,i_2,i_1)$ from $I$ such that 
\beq
{\hbox{
$i_k\ne i_{k+1}$ and $\sharp\{k: i_k=i\}=\ify$ for any $i\in I$.}}
\label{seq-con}
\eeq
Following Kashiwara \cite{K3}, 
we will associate to $\io$ a crystal structure
on $\ZZ^{\ify}$ and the embedding of crystals
\beq 
\Psi_{\io}:B(\ify)\hookrightarrow \ZZ^{\ify},
\label{psi}
\eeq
which we call the {\it Kashiwara embedding}.

The crystal structure 
on $\ZZ^{\ify}$ corresponding to $\io$ is defined as follows.
Let $\vec x = (\cd, x_k,\cd,x_2,x_1)\in 
\ZZ^{\ify}$.
For $k\geq1$, we set 
\beq
\sigma_k (\vec x):= x_k+\sum_{j>k}\lan h_{i_k},\al_{i_j}\ran x_j.
\label{sigma k}
\eeq
Since $x_j = 0$ for $j \gg 0$, the form $\sigma_k (\vec x)$ 
is well-defined, and $\sigma_k (\vec x) = 0 $ for $k \gg 0$.
For $i \in I$, let 
$\sigma^{(i)} (\vec x) := {\rm max}_{k: i_k = i} \sigma_k (\vec x)$,     and 
$$M^{(i)} = M^{(i)} (\vec x) := 
\{k: i_k = i, \sigma_k (\vec x) = \sigma^{(i)}(\vec x)\}.$$
Note that 
$\sigma^{(i)} (\vec x)\geq 0$, and that
$M^{(i)} = M^{(i)} (\vec x)$ is a finite set
if and only if $\sigma^{(i)} (\vec x) > 0$. 
Now we define the maps $\eit: \ZZ^{\ify} \lar \ZZ^{\ify} \sqcup\{0\}$
and $\fit: \ZZ^{\ify} \lar \ZZ^{\ify}$ by setting 
\beq
(\fit(\vec x))_k  = x_k + \delta_{k,{\rm min}\,M^{(i)}};
\label{action-f}
\eeq
\beq
(\eit(\vec x))_k  = x_k - \delta_{k,{\rm max}\,M^{(i)}} \,\, {\rm if}\,\,
\sigma^{(i)} (\vec x) > 0; \,\, {\rm otherwise} \,\, \eit(\vec x)=0.
\label{action-e}
\eeq

We also define the weight function and the functions 
$\vep_i$ and $\vp_i$ on $\ZZ^{\ify}$ by
$$wt(\vec x) := -\sum_{j=1}^{\ify} x_j \al_{i_j}, \,\,
\vep_i (\vec x) := \sigma^{(i)} (\vec x), \,\,
\vp_i (\vec x) := \lan h_i, wt(\vec x) \ran + \vep_i(\vec x).$$
An easy check shows that these maps make $\ZZ^{\ify}$
into a crystal.
We will denote this crystal by $\ZZ^{\ify}_{\io}$.
Note that, in general, the semigroup $\ZZ^{\ify}_{\geq 0}$ is not 
a subcrystal of $\ZZ^{\ify}_{\io}$ since it is not 
stable under the action of $\eit$'s.

The Kashiwara embedding is given by the following theorem.

\newtheorem{thm2}[df2]{Theorem}
\begin{thm2}
\label{emb}
There is a unique embedding of crystals 
\beq 
\Psi_{\io}:B(\ify)\hookrightarrow \ZZ^{\ify}_{\geq 0}
\subset \ZZ^{\ify}_{\io},
\label{psi}
\eeq
such that 
$\Psi_{\io} (u_{\ify}) = (\cd,0,\cd,0,0)$.
\end{thm2}

\nd
{\sl Proof.}\,\, 
The uniqueness of $\Psi_{\io}$ follows from the fact that every 
element of $B(\ify)$ is obtained from $u_{\ify}$ by a sequence of
operators $\fit$. 
To prove the existence, we show that $\Psi_{\io}$ 
can be obtained by iterating the following construction.
Recall that to every $i\in I$ we associate a crystal $B_i$ as in 
Example \ref{Example:crystal}.

\begin{thm2}[\cite{K3}]
\label{B=BBi}
For any $i\in I$, there is a unique embedding of crystals
\beqn
\Psi_i : B(\infty)&\hookrightarrow &B(\infty)\ot B_i,
\eeqn
such that $\Psi_i (u_{\ify}) = u_{\ify}\ot (0)_i$.
\end{thm2}

 An explicit formula for $\Psi_i$ is given as follows.
Let $x \mapsto x^*$ be the $\QQ(q)$-algebra antiautomorphism of 
$\uq$ given by:
\beq
e_i^*=e_i,\q f_i^*=f_i, \q (q^h)^*=q^{-h}.
\label{star}
\eeq
It is proved in \cite[Theorem 2.1.1]{K3} that this anti-automorphism
preserves $L(\ify)$, and that the induced action on $L(\ify)/qL(\ify)$
preserves the crystal $B(\ify)$.
Now for $b \in B(\ify)$ we have:
$\Psi_i(b)=b' \ot (-a)_i$,
where $a = \vep_i (b^*) \geq 0$ and $b' = (\eit^a (b^*))^*$.

Returning to the Kashiwara embedding $\Psi_{\io}$,
take any $b \in B(\ify)$ and define the elements $b_0, b_1, b_2, \cd$
of $B(\ify)$ and non-negative integers $a_1, a_2, \cd$ recursively by:
\beq
b_0 = b, \, \, \Psi_{i_k} (b_{k-1}) = b_k \ot (-a_k)_{i_k}
\q (k \geq 1).
\label{Kash recur}
\eeq
The definitions readily imply that $b_k = u_\ify$ and $a_k = 0$
for $k \gg 0$.
Thus, the sequence $(\cd,a_k,\cd,a_2,a_1)$ belongs to
$\ZZ^{\ify}_{\geq 0}$, and we set 
$$\Psi_{\io} (b) = (\cd,a_k,\cd,a_2,a_1).$$

The injectivity of $\Psi_{\io}$ follows from that of the $\Psi_{i_k}$.
To complete the proof of Theorem \ref{emb}, it remains to show that 
$\Psi_{\io}:B(\ify)\hookrightarrow \ZZ^{\ify}_{\io}$
is an embedding of crystals. 
Tracing the definition of $\Psi_{\io}$ and using Theorem \ref{B=BBi}, 
we only need to check the following:
if we identify a sequence 
$(\cd,0,0,a_k,\cd,a_2,a_1) \in \ZZ^{\ify}_{\io}$
with an element
$u_{\ify}\ot (-a_k)_{i_k} \ot \cd  \ot (-a_2)_{i_2}\ot (-a_1)_{i_1}$
of the tensor product of crystals 
$B(\ify)\ot B_{i_k} \ot \cd \ot B_{i_2} \ot B_{i_1}$
for some $k \gg 0$, then the crystal structure on 
$\ZZ^{\ify}_{\io}$ agrees with that of 
$B(\ify)\ot B_{i_k} \ot \cd \ot B_{i_2} \ot B_{i_1}$. 
This is a direct consequence of our definitions and 
Lemma 1.3.6 in \cite{K3}.
This concludes the proof of Theorem \ref{emb}. \qed

\vskip5pt

\nd
{\bf Remark.} \,\, The recursive definition (\ref{Kash recur})
of the Kashiwara embedding can be reformulated as follows:
if $\Psi_{\io} (b) = (\cd,a_k,\cd,a_2,a_1)$ then each $a_k$ is 
given by
\beq
a_k = \vep_{i_k} (\til e_{i_{k-1}}^{a_{k-1}} \cd 
\til e_{i_{1}}^{a_{1}} b^*) = 
{\rm max} \q \{a: \til e_{i_k}^{a} \til e_{i_{k-1}}^{a_{k-1}} \cd 
\til e_{i_{1}}^{a_{1}} b^* \neq 0\}.
\label{string analog}
\eeq
This is the crystal version of the {\it string parametrization}
introduced in Section~2 of \cite{BZ1}. 
It is not hard to show that, in the terminology 
of \cite{BZ1}, $\Psi_{\io} (b)$ is the {\it string} of $b^*$ 
(more precisely, of the global basis vector corresponding to
$b^*$) in direction $\io$. 
Passing from the $\eit$ to the $\fit$ transforms (\ref{string analog})
into one more equivalent description of the Kashiwara embedding:
the sequence $(\cd,a_k,\cd,a_2,a_1) = \Psi_{\io} (b)$
is uniquely determined by the conditions that 
\beq
b^* = \til f_{i_1}^{a_1} \til f_{i_{2}}^{a_{2}} \cd u_{\ify}, \,\,
{\rm and} \,\, 
\til e_{i_{k-1}}(\til f_{i_k}^{a_k} \til f_{i_{k+1}}^{a_{k+1}}
\cd u_{\ify})=0 \,\,{\rm for} \,\, k > 1.
\eeq

In the rest of the paper we deal with the following 

\vskip5pt

\nd
{\bf Main Problem.}\,\, Describe explicitly the image of $\Psi_{\io}$.

\section{Polyhedral Realizations of $B(\ify)$}
\setcounter{equation}{0}
\renewcommand{\theequation}{\thesection.\arabic{equation}}

\subsection{Piecewise-linear transformations $S_k$}
We will retain the notation introduced above.
In particular, we fix a sequence of indices
$\io:=(i_k)_{k\geq1}$ as in (\ref{seq-con}). 
Consider the infinite dimensional vector space 
$$
\QQ^{\ify}:=\{\vec{x}=
(\cd,x_k,\cd,x_2,x_1): x_k \in \QQ\,\,{\rm and }\,\,
x_k = 0\,\,{\rm for}\,\, k \gg 0\},
$$
and its dual $(\QQ^{\ify})^*:={\rm Hom}(\QQ^{\ify},\QQ)$. 
We will write a linear form $\vp \in (\QQ^{\ify})^*$ as
$\vp(\vec{x})=\sum_{k \geq 1} \vp_k x_k$ ($\vp_j\in \QQ$). 

For every $k \geq 1$, we define $\kp$ to be the minimal index $j$ 
such that $j > k$ and $i_j = i_k$. 
We also define $\km$ to be the maximal index $j$ 
such that $j < k$ and $i_j = i_k$; 
if $i_j \neq i_k$ for $1 \leq j < k$, then we set $\km = 0$.
Let $\beta_k \in (\QQ^{\ify})^*$ be a linear form given by
\beq
\beta_k (\vec x) = \sigma_k (\vec x) - \sigma_{\kp} (\vec x),
\label{beta}
\eeq
where the forms $\sigma_k$ are defined by (\ref{sigma k}). 
Since $\lan h_{i},\al_{i}\ran   = 2$ for any $i \in I$, we have
\beq
\beta_k(\vec{ x}) =   x_k+\sum_{k<j<\kp}
\lan h_{i_k},\al_{i_j}\ran x_j+x_{\kp}.
\label{beta+}
\eeq
We will also use the convention that $\beta_0 (\vec x) = 0$
for all $\vec x \in \QQ^{\ify}$. 
Using this notation, for every $k \geq 1$, we define a 
piecewise-linear operator 
$S_k = S_{k,\io}$ on ${(\QQ^{\ify})}^*$ by
\beq 
S_k\,(\vp) :=\cases{
\vp - \vp_k \beta_k & if $\vp_k > 0$,\cr
\vp - \vp_k \beta_{\km} & if $\vp_k \leq 0$.\cr}
\label{S_k}
\eeq 
An easy check shows that 
$(S_k)^2=S_k.$

\subsection{Main theorem}
For a sequence $\io=(i_k)_{k\geq1}$ satisfying (\ref{seq-con}),
we denote by $\Xi_{\io} \subset {(\QQ^{\ify})}^*$
the subset of linear forms that are obtained from the coordinate forms 
$x_j$ by applying transformations $S_k = S_{k,\io}$ (see (\ref{S_k})).
In other words, we set
\beq
\Xi_{\io} :=  
\{S_{j_l}\cd S_{j_2}S_{j_1}x_{j_0}: l\geq 0,\,\, j_0, \cd, j_l \geq 1\}.
\label{Xi}
\eeq

\nd
Recall that, for $k \geq 1$,
the condition $\km = 0$ means that $i_j \neq i_k$ for $1 \leq j < k$.
We will impose on $\io$ the following positivity assumption: 
\beq
{\rm if} \,\, \km = 0 \,\, {\rm then}\,\, \vp_k \geq 0 \,\,
{\rm for \,\, any} \,\, \vp = \sum\vp_j x_j\in \Xi_{\io}. 
\label{positivity}
\eeq

Now we are in a position to formulate our main result.

\newtheorem{thm3}{Theorem}[section]
\begin{thm3}
\label{main}
Let $\io$ be a sequence of indices satisfying (\ref{seq-con})
and the positivity assumption (\ref{positivity}). 
Let 
$\Psi_{\io}:B(\ify)\hookrightarrow \ZZ^{\ify}_{\geq 0}$
be the corresponding Kashiwara embedding.
Then the image ${\rm Im} \,(\Psi_{\io})$ is equal to 
\beq
\Sigma_{\io} := \{x\in\ZZ^{\ify}_{\geq0}\subset \QQ^{\ify}:
\vp(x)\geq0 \,\,{\rm for\,\, any }\,\,\vp\in \Xi_{\io}\}.
\label{Sigma}
\eeq
\end{thm3}  
 
\vskip5pt
\nd
{\sl Proof.}\,\,
In view of Theorem \ref{emb}, the image ${\rm Im} \,(\Psi_{\io})$
is a subcrystal of $\ZZ^{\ify}_{\io}$ obtained 
by applying the operators $\fit$ to 
$\Psi_{\io}(u_{\ify})=\vec{0}=(\cd,0,0,0)$;
in particular, 
${\rm Im} \,(\Psi_{\io}) \subset \ZZ^{\ify}_{\geq 0}$.
Since $\vec{0} \in \Sigma_{\io}$, 
the inclusion ${\rm Im} \,(\Psi_{\io}) \subset \Sigma_{\io}$
follows from the fact that $\Sigma_{\io}$ is closed under all $\fit$. 
Let us prove this fact.
Let $\vec x = (\cd,x_2,x_1)\in \Sigma_{\io}$ and $i\in I$,
and suppose that $\fit\vec x = (\cd,x_k+1,\cd,x_2,x_1)$
(in particular, $i_k=i$). 
We need to show that
\beq
\vp(\fit\vec x) \geq0
\label{vpfit}
\eeq
for any $\vp = \sum\vp_jx_j\in\Xi_{\io}$. 
Since $\vp(\fit\vec x)=\vp(\vec x)+\vp_k \geq \vp_k$, 
it is enough to consider the case when $\vp_k < 0$. 
By (\ref{positivity}), we have $\km \geq 1$.
Remembering (\ref{action-f}), we have 
$\sigma_k(\vec x)>\sigma_{\km}(\vec x)$ and then by 
(\ref{beta}), we conclude that 
$$
\beta_{\km} (\vec x)= 
\sigma_{\km}(\vec x) - \sigma_{k}(\vec x)\leq -1.
$$
It follows that 
\beqn
\vp(\fit\vec x) & = & \vp(\vec x)+\vp_k\nn\\
		  & \geq & \vp(\vec x)-\vp_k\beta_{\km}(\vec x)\nn\\
		  & = & (S_k\vp)(\vec x) \geq 0,\nn
\eeqn
since $S_k\vp \in \Xi_{\io}$. 
This proves the inclusion ${\rm Im} \,(\Psi_{\io}) \subset \Sigma_{\io}$.

To prove the reverse inclusion 
$\Sigma_{\io} \subset {\rm Im} \,(\Psi_{\io})$,
we first show that $\Sigma_{\io}$ is a subcrystal of $\ZZ^{\ify}_{\io}$,
i.e., that $\eit \Sigma_{\io}\subset \Sigma_{\io}\sqcup\{0\}$ 
for any $i\in I$.
Let $\vec x = (\cd,x_2,x_1)\in \Sigma_{\io}$ and $i\in I$, and  
suppose that 
$\eit\vec x = (\cd,x_k-1,\cd,x_2,x_1)$; in particular, $i_k=i$. 
We need to show that
\beq
\vp(\eit\vec x) \geq 0,
\label{vpeit}
\eeq
for any $\vp = \sum\vp_jx_j\in\Xi_{\io}$.
Since $\vp(\eit\vec x) = \vp(\vec x)-\vp_k \geq -\vp_k$, 
it is enough to consider the case when $\vp_k > 0$. 
Remembering (\ref{action-e}), 
we have $\sigma_k(\vec x)>\sigma_{\kp}(\vec x)$  and then by
(\ref{beta}), we conclude that 
$$
\beta_k(\vec x)
= \sigma_{k}(\vec x) - \sigma_{\kp}(\vec x) \geq 1.
$$
It follows that 
\beqn
\vp(\eit\vec x) & = & \vp(\vec x)-\vp_k\nn\\
		  & \geq & \vp(\vec x)-\vp_k\beta_k (\vec x)\nn\\
		  & = & (S_k\vp)(\vec x) \geq 0,\nn
\eeqn
since $S_k\vp \in \Xi_{\io}$. 

To complete the proof of the inclusion 
$\Sigma_{\io} \subset {\rm Im} \,(\Psi_{\io})$,
we make the following observation:
if $\vec x \in \ZZ^{\ify}_{\geq 0, \io}$ and $\vec x \neq \vec 0$
then $\eit \vec x \neq 0$ for some $i \in I$.
(Indeed, one can take $i = i_j$, where $j$ is the maximal index such that
$x_j > 0$.).
Since $\Sigma_{\io} \subset \ZZ^{\ify}_{\geq 0}$, we conclude that 
every $\vec x \in \Sigma_{\io}$ can be transformed to $\vec 0$
by a sequence of operators $\eit$. 
By (\ref{eeff}), $\vec x$ is obtained from $\vec 0$
by a sequence of operators $\fit$, hence belongs to 
${\rm Im} \,(\Psi_{\io})$, and we are done. \qed

\subsection{Remarks on the positivity assumption}
In this subsection, we shall give some equivalent reformulations
of the positivity assumption (\ref{positivity}).
This will allow us to sharpen Theorem \ref{main}.
We retain all the previous notation; in particular,
we fix a sequence $\io$ satisfying (\ref{seq-con}).

For every $k \geq 1$, we introduce the transformations
$E_k, F_k: \ZZ^{\ify}\longrightarrow \ZZ^{\ify} \sqcup\{0\}$ 
that act on $\vec x=(\cd,x_2,x_1)$ by
\beqn
E_k(\vec x)& := & \cases{
(\cd, x_k-1,\cd, x_2, x_1) & if $\beta_k(\vec x)>0$,\cr
\q 0 & otherwise.\cr}
\label{E}\\
F_k(\vec x)& := & \cases{
(\cd, x_k+1,\cd,x_2, x_1) & if $\beta_{\km} (\vec x)<0$ or $\km = 0$,\cr
\q 0 & otherwise.\cr}
\label{F}
\eeqn
Comparing these definitions with (\ref{action-f}) and (\ref{action-e}),
we see that the operators $\eit$ and $\fit$ can be written as
\beq
\eit = E_{{\rm max}\,M^{(i)}}, \,\, \fit = F_{{\rm min}\,M^{(i)}}.
\label{local form}
\eeq
Let $\Phi = \Phi_\io \subset \ZZ^{\ify}$ be the set 
of all $\vec x$ obtained from $\vec 0 = (\cd, 0,0)$ by 
applying transformations $E_k$ and $F_k$.
Let $\Phi^+ \subset \Phi$ be the set of all $\vec x$ obtained from 
$\vec 0$ by applying the $F_k$.
Recall also the definitions (\ref{Xi}) of the set of linear forms
$\Xi = \Xi_\io$, and (\ref{Sigma}) of the subset 
$\Sigma = \Sigma_\io \subset \ZZ^{\ify}_{\geq0}$.
We will say that a transformation $S_k$ acts positively on a linear form
$\vp$ if $\vp_k > 0$, i.e., the first possibility in (\ref{S_k})
is realized. 
We define $\Xi^+ \subset \Xi$ to be the set of forms 
obtained from the coordinate forms $x_j$ by applying positive actions
of the transformations $S_k$. 
Let 
\beq
\Sigma^+ = \{x \in \ZZ^{\ify}_{\geq0}: 
\vp(x)\geq0 \,\,{\rm for\,\, any }\,\,\vp\in \Xi^+\}.
\label{Sigma^+}
\eeq
The sets $\Phi, \Phi^+, \Sigma$, and $\Sigma^+$ are related to each other
and to the image ${\rm Im} \,(\Psi_{\io})$ of the Kashiwara embedding
as follows (note that we do not assume (\ref{positivity}) here). 

\newtheorem{pr3}[thm3]{Proposition}
\begin{pr3}
We have
\beq					
\Sigma \subset \Sigma^+ \subset {\rm Im} \,(\Psi_{\io}) \subset
\Phi^+ \subset \Phi.
\label{5 inclusions}
\eeq
\end{pr3}
{\sl Proof.\,\,}
The inclusions $\Sigma \subset \Sigma^+$ and  $\Phi^+ \subset \Phi$
are obvious.
Since every $\vec x \in {\rm Im} \,(\Psi_{\io})$ is obtained
from $\vec 0$ by applying the $\fit$, the inclusion 
${\rm Im} \,(\Psi_{\io}) \subset \Phi^+$ follows from 
the second equality in (\ref{local form}). 
The remaining inclusion $\Sigma^+ \subset {\rm Im} \,(\Psi_{\io})$
is proved by the same argument as the inclusion 
$\Sigma_{\io} \subset {\rm Im} \,(\Psi_{\io})$ in Theorem \ref{main}
(it is seen by inspection that this argument does not use 
(\ref{positivity})). \qed

\vskip5pt
\nd
{\bf Remark.} \,\, Note that the same argument as in the proof of 
Theorem \ref{main} establishes the following inclusions:
\beq
E_k \Sigma \subset \Sigma \sqcup \{0\}, \,\, 
E_k \Sigma^+ \subset \Sigma^+ \sqcup \{0\} \,\, {\rm for \,\, any} \,\, 
k \geq 1.
\label{E-inv}
\eeq

\vskip5pt

\begin{thm3}
\label{posit char}
The following conditions are equivalent:
\begin{enumerate}
\item
All the inclusions in (\ref{5 inclusions}) are equalities.
\item
$\Phi \subset \ZZ^{\ify}_{\geq0}$.
\item
$\Xi$ satisfies the positivity assumption (\ref{positivity}).
\end{enumerate}
\end{thm3}
{\sl Proof.\,\,}
The implication (i) $\Rightarrow$ (ii) is obvious
since $\Phi^+ \subset \ZZ^{\ify}_{\geq0}$.
To prove (iii) $\Rightarrow$ (i), we notice that 
(\ref{positivity}) implies the following companion of (\ref{E-inv}):
\beq
F_k \Sigma \subset \Sigma \sqcup \{0\} \,\, {\rm for \,\, any} \,\, 
k \geq 1;
\label{F-inv}
\eeq
this is proved by the same argument as the fact that 
$\Sigma$ is closed under all $\fit$ in the proof of Theorem \ref{main}. 
Combining (\ref{E-inv}) and (\ref{F-inv}), we conclude that
$\Phi \subset \Sigma$, which implies (i). 

It remains to prove (ii) $\Rightarrow$ (iii). 
We will deduce this from the following lemma.

\newtheorem{lm3}[thm3]{Lemma}
\begin{lm3}
\label{S_k inv}
The set of linear forms that take nonnegative values on $\Phi$,
is closed under all transformations $S_k$. 
\end{lm3}
{\sl Proof.\,\,}
Suppose $\vp (\vec x) = \sum_{k \geq 1} \vp_k x_k \geq 0$
for all $\vec x \in \Phi$. 
Take any $\vec x \in \Phi$.
We need to show that $S_k \vp (\vec x) \geq 0$ for any $k \geq 1$.
First suppose that $\vp_k > 0$, i.e., $S_k \vp = \vp - \vp_k \beta_k$. 
If $\beta_k (\vec x) \leq 0$ then 
$S_k \vp (\vec x) \geq \vp (\vec x) \geq 0$;
so we can assume that $\beta_k (\vec x) = l > 0$. 
Using (\ref{beta+}) and (\ref{E}), we conclude that 
$(E_k)^l \vec x \in \Phi$, and
$$S_k \vp (\vec x) = \vp (\vec x) - l \vp_k = \vp ((E_k)^l \vec x) \geq 0,$$
as required.

It remains to consider the case when $\vp_k < 0$, i.e., 
$S_k \vp = \vp - \vp_k \beta_{\km}$. 
If $\beta_k (\vec x) \geq 0$ then 
$S_k \vp (\vec x) \geq \vp (\vec x) \geq 0$;
so we can assume that $\beta_k (\vec x) = -l < 0$. 
Using (\ref{beta+}) and (\ref{F}), we conclude that 
$(F_k)^l \vec x \in \Phi$, and
$$S_k \vp (\vec x) = \vp (\vec x) + l \vp_k = \vp ((F_k)^l \vec x) \geq 0,$$
as required. \qed

\vskip 5pt

Now we can complete the proof of (ii) $\Rightarrow$ (iii). 
By (ii) and Lemma \ref{S_k inv}, every form $\vp \in \Xi$
takes nonnegative values on $\Phi$
since $\vp$	is obtained from some coordinate form $x_j$ by applying 
the transformations $S_k$, and $x_j$ is nonnegative on $\Phi$.
In particular, if $\km = 0$ for some $k \geq 1$ then
$$\vp_k = \vp (F_k \vec 0) \geq 0,$$
which proves (\ref{positivity}).
Theorem \ref{posit char} is proved. \qed

\vskip5pt

\nd
{\bf Remarks.} (a) Using Theorem \ref{posit char}, one can produce 
several other equivalent reformulations of (\ref{positivity}).
For instance, each of the following two conditions is also
equivalent to (\ref{positivity}):

\nd
(iv) $\Sigma = \Sigma^+$;

\vskip5pt

\nd
(v) $\Phi = \Phi^+$

\vskip5pt

\nd
(the implications (i) $\Rightarrow$ (iv) $\Rightarrow$ (iii)
and (i) $\Rightarrow$ (v) $\Rightarrow$ (ii) are obvious). 

\vskip5pt

(b) It would be interesting to know if (\ref{positivity})
holds for any symmetrizable Kac-Moody algebra and any sequence $\io$. 
This will be true in all the examples considered in the rest of the paper. 
In fact, in all these examples we will have $\Xi = \Xi^+$,
which is stronger than the condition (iv) above. 

\subsection{Periodic case}
In the subsequent sections we shall only treat the following special 
infinite sequence $\io$. 
We fix some linear ordering of the index set $I$, i.e., 
identify $I$ with $\{1,2,\cd,n\}$.
Then we take
$$
\io = (\cd,\underbrace{n,\cd,2,1}_{},
\cd,\underbrace{n,\cd,2,1}_{},\underbrace{n,\cd,2,1}_{}).
$$
In other words, $i_k = \overline k$, where 
$\overline k \in \{1,2,\cd,n\}$ is congruent to $k$ modulo $n$.
We call this sequence $\io$ {\it periodic}.
Relative to the periodic sequence, the above notation 
simplifies as follows.

\nd
First of all, for any $k\geq 1$ we have $\kp=k+n$;
we also have $\km=k-n$ if $k>n$, and $\km=0$ if $k\leq n$.
The forms $\beta_k$ take the form
$$
\beta_k(\vec x)=x_k+\sum_{j=k+1}^{k+n-1}
\lan h_{\overline k},\al_{\overline j}\ran x_j+x_{k+n},
$$
and the transformations $S_k$ can be written as
\beq 
S_k\,(\vp) :=\cases{
\vp - \vp_k \beta_k & if $\vp_k \geq 0$,\cr
\vp - \vp_k\beta_{k-n} & if $k > n, \, \vp_k <  0$,\cr
\vp     & if $1 \leq k \leq n, \, \vp_k < 0$.\cr}
\label{S_k-periodic}
\eeq 
Finally, the positivity assumption (\ref{positivity}) in Theorem \ref{main}
takes the form 
\beq
\vp_i \geq 0 \,\,{\rm for \,\, any} \,\, i = 1,2,\cd,n \,\,{\rm and} \,\,
\vp = \sum\vp_j x_j\in \Xi_{\io}. 
\label{perpos}
\eeq
This means that, for $\vp \in \Xi_{\io}$, the third opportunity 
in (\ref{S_k-periodic}) is never realized.

\section{Rank 2 case}
\setcounter{equation}{0}
\renewcommand{\theequation}{\thesection.\arabic{equation}}

In this section, we specialize Theorem \ref{main} to the Kac-Moody
algebras of rank 2. 
We will give an explicit description of the image of the Kashiwara 
embedding. 
This description sharpens the one given by
Kashiwara in \cite[Sect.2]{K3}.

Without loss of generality, we can and will assume that $I=\{1,2\}$,
and $\io = (\cd,2,1,2,1)$. 
Let the Cartan data be given by:
$$
\lan h_1,\al_1\ran= \lan h_2,\al_2\ran=2, \,\, \lan h_1,\al_2\ran=-c_1, 
\,\, \lan h_2,\al_1\ran=-c_2.
$$
Here we either have $c_1 = c_2 = 0$, or both $c_1$ and $c_2$ are 
positive integers. 
We set $\lambda = c_1 c_2 - 2$, and define the integer sequence 
$a_l = a_l (c_1, c_2)$ for $l \geq 0$ by setting $a_0 = 0, \, a_1 = 1$
and, for $k \geq 1$,
\beq
a_{2k}  = c_1 P_{k-1} (\lambda), \,\, 
a_{2k+1} = P_k (\lambda) + P_{k-1} (\lambda),
\label{defcoeff}
\eeq
where the $P_k (\lambda)$ are {\it Chebyshev polynomials} given by
\beq
 P_k (\al + \al^{-1}) = 
 {\al^{k+1} - \alpha^{-k-1}     \over \alpha - \alpha^{-1}}.
\label{cheb}
\eeq
Equivalently, the generating function for Chebyshev polynomials is given by
\beq
 \sum_{k \geq 0} P_k (\lambda) z^k = (1 - \lambda z + z^2)^{-1}.
\label{gen-cheb}
\eeq
The several first Chebyshev polynomials and terms $a_l$ are given by
$$P_0 (\lambda) = 1, \, P_1 (\lambda) = \lambda, \,
P_2 (\lambda) = \lambda^2 - 1, \, 
P_3 (\lambda) = \lambda^3 - 2 \lambda,$$
$$a_2 = c_1, \, a_3 = c_1 c_2 - 1, \, a_4 = c_1 (c_1 c_2 - 2),$$
$$a_5 = (c_1 c_2 - 1)(c_1 c_2 - 2) - 1, \, 
a_6 = c_1 (c_1 c_2 - 1)(c_1 c_2 - 3),$$
$$a_7 = c_1 c_2 (c_1 c_2 - 2)(c_1 c_2 - 3) - 1.$$
Let $l_{\rm max} = l_{\rm max} (c_1, c_2)$ be the minimal index
$l$ such that $a_{l+1} < 0$
(if $a_l \geq 0$ for all $l \geq 0$,
then we set $l_{\rm max} = + \infty$). 
By inspection, if $c_1 c_2 = 0$ (resp. $1,2,3$) 
then $l_{\rm max} = 2$ (resp. $3, 4, 6$).
Furthermore, if $c_1 c_2 \leq 3$ then $a_{l_{\rm max}} = 0$ and 
$a_l > 0$ for $1 \leq l < l_{\rm max}$.
On the other hand, if $c_1 c_2 \geq 4$, i.e., $\lambda \geq 2$,
it is easy to see from (\ref{cheb}) or (\ref{gen-cheb})
that $P_k (\lambda) > 0$ for $k \geq 0$, hence 
$a_l > 0$ for $l \geq 1$; in particular, in this case
$l_{\rm max} = + \infty$.  

\newtheorem{thm4}{Theorem}[section]
\begin{thm4}
\label{rank 2-thm}
In the rank 2 case, the image of the Kashiwara embedding
is given by
\beq
{\rm Im} \,(\Psi_{\io}) = \left\{(\cd,x_2,x_1)\in\ZZ_{\geq0}^{\ify}: 
\begin{array}{l}
x_k = 0 \,\,{\rm for}\,\,k > l_{\rm max}, \\
a_l x_l -a_{l-1} x_{l+1} \geq0 \,\,{\rm for} \,\,
1 \leq l < l_{\rm max}
\end{array}
\right\}.
\label{rank 2-cone}
\eeq
\end{thm4} 
\vskip5pt
\nd
{\sl Proof.}\,\,
We will deduce our theorem from Theorem \ref{main}.
Thus, our first goal is to describe the set of linear forms 
$\Xi_{\io}$     (see (\ref{Xi})), and to check 
the positivity assumption (\ref{perpos}). 
For $k \geq 1$ and $0 \leq l < l_{\rm max}$, we set
\beq
\vp^{(l)}_k = S_{k+l-1}\cd S_{k+1}S_{k} x_k;
\label{S-forms}
\eeq
in particular, $\vp^{(0)}_k = x_k$. 
We also define $a'_l = a_l (c_2, c_1)$, i.e., the numbers $a'_l$ 
are given by (\ref{defcoeff}) with $c_1$ replaced by $c_2$.

\newtheorem{lm4}[thm4]{Lemma}
\begin{lm4}
\label{rank2-lm1}
\begin{enumerate}
\item [{\rm (i)}] If $k$ is odd then 
$\vp^{(l)}_k = a_{l+1} x_{k+l} - a_l x_{k+l+1}$; 
if $k$ is even then 
$\vp^{(l)}_k = a'_{l+1} x_{k+l} - a'_l x_{k+l+1}$.
\item [{\rm (ii)}] If $c_1 c_2 \leq 3$, i.e., $l_{\rm max} < + \infty$,
then $\vp^{(l_{\rm max}-1)}_k  = -x_{k+l_{\rm max}}$.
\item [{\rm (iii)}] The set $\Xi_{\io}$ consists of all linear forms
$\vp^{(l)}_k$ with $k \geq 1$ and $0 \leq l < l_{\rm max}$.
\item [{\rm (iv)}] The set $\Xi_{\io}$ satisfies the positivity assumption 
(\ref{perpos}).
\end{enumerate}
\end{lm4}

\nd
{\sl Proof.\,\,}  (i) In view of periodicity, it is enough to show that 
$\vp^{(l)}_1 = a_{l+1} x_{l+1} - a_l x_{l+2}$ for $0 \leq l < l_{\rm max}$.
We prove this by induction on $l$. 
The claim is obviously true for $l=0$, since $a_0 = 0$ and $a_1 = 1$.
So let us assume that the claim  is true for some $\vp^{(l)}_1$ such that
both $a_l$ and $a_{l+1}$ are positive; we need to show that
the claim is then true for $\vp^{(l+1)}_1$.      
By (\ref{S_k-periodic}),
$$\vp^{(l+1)}_1 = S_{l+1} \vp^{(l)}_1 =  
\vp^{(l)}_1 - a_{l+1} \beta_{l+1} = 
a_{l+1} x_{l+1} - a_l x_{l+2} - a_{l+1} \beta_{l+1},$$
where the forms $\beta_l$ are given by
\beqnn
\beta_{2k-1}(\vec x)& = & x_{2k-1}-c_1 x_{2k}+x_{2k+1},\\
\beta_{2k}(\vec x)& = & x_{2k}-c_2x_{2k+1}+x_{2k+2}.
\eeqnn
Therefore, $\vp^{(l+1)}_1 = (c_1 a_{l+1} - a_l) x_{l+2} - a_{l+1} x_{l+3}$
if $l$ is even, and 
$\vp^{(l+1)}_1 = (c_2 a_{l+1} - a_l) x_{l+2} - a_{l+1} x_{l+3}$
if $l$ is odd. 
It remains to show that the sequence $(a_l)$ satisfies the recursions
$$a_{2k+2}=c_1 a_{2k+1}-a_{2k}, \,\, a_{2k+1}=c_2 a_{2k}-a_{2k-1}.$$
These recursions follow from (\ref{defcoeff})
with the help of the well-known recursion 
$P_k (\lambda) = \lambda P_{k-1} - P_{k-2}$
for Chebyshev polynomials (the latter recursion is an easy consequence
of (\ref{gen-cheb})). 

(ii) Again it is enough to treat the case $k = 1$. 
By part (i), we have 
$\vp^{(l_{\rm max}-1)}_1  = a_{l_{\rm max}} x_{l_{\rm max}} - 
a_{l_{\rm max}-1} x_{l_{\rm max}+1}$.
So we only need to show that $a_{l_{\rm max}} = 0$ 
and $a_{l_{\rm max}-1} = 1$ in all the cases when $c_1 c_2 \leq 3$.
This is just seen by inspection.  

(iii) We only need to show that the set of all linear forms
$\vp^{(l)}_k$ with $k \geq 1$ and $0 \leq l < l_{\rm max}$
is closed under all the transformations $S_j$. 
If $l=0$ then the only $S_j$ that acts non-trivially on 
$\vp^{(0)}_k = x_k$     is $S_k$, and we have $S_k x_k = \vp^{(1)}_k$.
If $1 \leq l < l_{\rm max}$ then, in view of (i), only
$S_{k+l}$ and $S_{k+l+1}$ can act non-trivially on $\vp^{(l)}_k$.
By definition, $S_{k+l} \vp^{(l)}_k = \vp^{(l+1)}_k$.
We complete the proof by showing that 
$S_{k+l+1} \vp^{(l)}_k = \vp^{(l-1)}_k$.
Once again, using periodicity we can assume that $k=1$. 
Using (i) and (\ref{S_k-periodic}), we have
$$S_{l+2} \vp^{(l)}_1 = \vp^{(l)}_1 + a_l \beta_l = 
S_l \vp^{(l-1)}_1 + a_l \beta_l = (\vp^{(l-1)}_1 - a_l \beta_l) 
+ a_l \beta_l = \vp^{(l-1)}_1,$$
as claimed. 

Finally, part (iv) is an immediate consequence of (i) and (iii). \qed

\vskip5pt

Now we return to the proof of Theorem \ref{rank 2-thm}.
Using Theorem \ref{main} and parts (iii) and (iv) of Lemma \ref{rank2-lm1},
we conclude that
\beq
{\rm Im} \,(\Psi_{\io}) = \Sigma_\io = 
\{(\cd,x_2,x_1)\in\ZZ_{\geq0}^{\ify}: 
\vp^{(l)}_k \geq 0 \,\,{\rm for}\,\, k \geq 1, \, 0 \leq l < l_{\rm max}\}.
\label{intermed}
\eeq
Comparing this with the desired answer (\ref{rank 2-cone}), and using 
parts (i) and (ii) of Lemma \ref{rank2-lm1},
it only remains to show that the inequalities $\vp^{(l)}_k \geq 0$
in (\ref{intermed}) are redundant when $k > 1$ and $l < l_{\rm max}-1$,
that is, they are consequences of the remaining inequalities.
We will prove this by showing that the inequality $\vp^{(l-1)}_{k} \geq 0$ 
with $k > 1$ and $l < l_{\rm max}$ is a consequence of 
$\vp^{(l)}_{k-1} \geq 0$. 
As above, by using periodicity, it suffices to show that 
$\vp^{(l)}_{1} \geq 0$ implies $\vp^{(l-1)}_{2} \geq 0$. 
By Lemma \ref{rank2-lm1} (i), we have
$$\vp^{(l)}_1 = a_{l+1} x_{l+1} - a_l x_{l+2}, \,\,
\vp^{(l-1)}_2 = a'_{l} x_{l+1} - a'_{l-1} x_{l+2},$$
which easily implies that
$$a_{l+1} \vp^{(l-1)}_2 = a'_l \vp^{(l)}_1 + 
(a_l a'_l - a_{l+1}a'_{l-1}) x_{l+2}.$$
To complete the proof, it suffices to show that 
$a_l a'_l - a_{l+1}a'_{l-1}     > 0$ for all $l \geq 1$.
In fact, we claim that the numbers $a_l$ and $a'_l$
satisfy the identity
\beq
a_l a'_l - a_{l+1}a'_{l-1} = 1
\label{ident}
\eeq
for $l \geq 1$; this is a consequence of (\ref{defcoeff}) 
and the following identities for Chebyshev polynomials:
$$(\lambda + 2) P_k (\lambda)^2 - 
(P_{k+1} (\lambda) + P_{k} (\lambda)) (P_{k} (\lambda) + P_{k-1} (\lambda))
= 1,$$
$$(P_{k} (\lambda) + P_{k-1} (\lambda))^2 - 
(\lambda + 2) P_{k} (\lambda) P_{k-1} (\lambda) = 1$$
(the latter identities follow readily from (\ref{cheb})). 
Theorem \ref{rank 2-thm} is proved. \qed

\vskip5pt

Note that the cases when $l_{\rm max} < +\infty$, i.e., when the
image ${\rm Im} \,(\Psi_{\io})$ is contained in a lattice of finite rank,
are precisely those when the Kac-Moody algebra $\ge$ is of finite type.
If $\ge$ is of type $A_1 \times A_1$ (resp. $A_2$, $B_2$ or $C_2$, $G_2$)
then $l_{\rm max} = 2$ (resp. $3, 4, 6$).
Not surprisingly, in each case $l_{\rm max}$ is the number of positive
roots of $\ge$. 

In conclusion of this section, we illustrate Theorem \ref{rank 2-thm}
by the example when $c_1 = c_2 = 2$, i.e., $\ge$ is the affine Kac-Moody
Lie algebra of type $A^{(1)}_1$. 
In this case, we have $\lambda = c_1 c_2 - 2 = 2$.
It follows at once from (\ref{gen-cheb}) that $P_k (2) = k+1$;
hence, (\ref{defcoeff}) gives $a_l = l$ for $l \geq 0$.
We see that for type $A^{(1)}_1$, the  image of the Kashiwara embedding
is given by
\beq
{\rm Im} \,(\Psi_{\io}) = \{(\cd,x_2,x_1)\in\ZZ_{\geq0}^{\ify}: 
l x_l - (l-1) x_{l+1} \geq0 \,\,{\rm for}\,\, l \geq 1\}.
\label{A^{(1)}_1}
\eeq

\section{$A_n$-case}
\setcounter{equation}{0}
\renewcommand{\theequation}{\thesection.\arabic{equation}}
In this section we shall apply 
Theorem \ref{main} to the case when $\ge= \ssl_{n+1}$ is of type $A_n$. 
We will identify the index set $I$ with $[1,n] := \{1,2,\cd,n\}$ 
in the standard way; thus, the Cartan matrix 
$(a_{i,j}= \lan h_i,\al_j\ran )_{1 \leq i,j \leq n}$ is given by 
$a_{i,i}=2$, $a_{i,j}=-1$ for $|i-j|=1$, and 
$a_{i,j}=0$ otherwise. 
We will find the image of the Kashiwara embedding 
${\rm Im} \,(\Psi_{\io}) \subset \ZZ^{\ify}_{\geq 0}$
for the periodic sequence 
$$
\io = (\cd,\underbrace{n,\cd,2,1}_{},
\cd,\underbrace{n,\cd,2,1}_{},\underbrace{n,\cd,2,1}_{}).
$$

To formulate the answer, it will be convenient for us 
to change the indexing set for $\ZZ^{\ify}_{\geq 0}$
from $\ZZ_{\geq 1}$ to $\ZZ_{\geq 1} \times [1,n]$.
We will do this with the help of the bijection 
$\ZZ_{\geq 1} \times [1,n] \to \ZZ_{\geq 1}$ given by
$(j;i) \mapsto (j-1)n + i$. 
Thus, we will write an element $\vec x \in \ZZ^{\ify}_{\geq 0}$
as a doubly-indexed family $(x_{j;i})_{j \geq 1, i \in [1,n]}$
of nonnegative integers.
We will adopt the convention that $x_{j;i} = 0$ unless
$j \geq 1$ and  $i \in [1,n]$; in particular, $x_{j;0} = x_{j;n+1} = 0$
for all $j$. 

\newtheorem{thm5}{Theorem}[section]
\begin{thm5} 
\label{A_n}
In the above notation, the image ${\rm Im} \,(\Psi_{\io})$ 
of the Kashiwara embedding is the set of all integer families
$(x_{j;i})$ such that $x_{j;i} = 0$ for $i+j > n+1$, and 
$x_{1;i} \geq x_{2;i-1} \geq \cd \geq x_{i;1} \geq 0$ 
for $1 \leq i \leq n$. 
\end{thm5}
\vskip5pt
\nd
{\sl Proof.}\,\,
We will follow the proof of Theorem \ref{rank 2-thm}.
So we first describe the set of linear forms 
$\Xi_{\io}$, and check the positivity assumption (\ref{perpos}). 
To do this, we need some terminology and notation. 
For $(j;i) \in \ZZ_{\geq 1} \times [1,n]$, we will write 
the piecewise-linear transformation $S_{(j-1)n+i}$ 
as $S_{j;i}$; if $(j;i) \notin \ZZ_{\geq 1} \times [1,n]$ 
then $S_{j;i}$ is understood as the identity transformation. 
For $l \geq 0$ we set
\beq
S^{(l)}_{j;i} := S_{j;i+l-1}\cd S_{j;i+1} S_{j;i}.
\label{S-strings}
\eeq
(again with the understanding that $S^{(0)}_{j;i}$ is 
the identity transformation). 
For $i \in [1,n]$, by an $i$-{\it admissible partition}
we will mean an integer sequence $\lambda = (\lambda_1, \cd, \lambda_i)$
such that $n+1-i \geq \lambda_1 \geq \cd \geq \lambda_i \geq 0$
(if we represent partitions by Young diagrams in a usual way,
then this condition means that the diagram of $\lambda$ fits into
the $i \times (n+1-i)$ rectangle). 
For every $(j;i) \in \ZZ_{\geq 1} \times [1,n]$ and an 
$i$-admissible partition $\lambda$, we define the linear form 
$\vp^{(\lambda)}_{j;i}$ by
\beq
\vp^{(\lambda)}_{j;i} = 
S^{(\lambda_i)}_{j+i-1;1}\cd S^{(\lambda_2)}_{j+1;i-1}
S^{(\lambda_1)}_{j;i} x_{j;i}.
\label{S-forms A_n}
\eeq

\newtheorem{lm5}[thm5]{Lemma}
\begin{lm5}
\label{lm5.1}
\begin{enumerate}
\item [{\rm (i)}] The forms $\vp^{(\lambda)}_{j;i}$ are given by
\beq
\vp^{(\lambda)}_{j;i} = 
\sum_{k=1}^i (x_{j+k-1;i-k+1+\lambda_k} - x_{j+k;i-k+\lambda_k}).
\label{A_n exp}
\eeq
\item [{\rm (ii)}] If $\lambda_k = n+1-i$ for $k = 1, \cd, i$
then $\vp^{(\lambda)}_{j;i} = -x_{j+i;n+1-i}$.
\item [{\rm (iii)}] The set $\Xi_{\io}$ consists of all linear forms
$\vp^{(\lambda)}_{j;i}$, where 
$(j;i) \in \ZZ_{\geq 1} \times [1,n]$ and $\lambda$ is an 
$i$-admissible partition. 
\item [{\rm (iv)}] The set $\Xi_{\io}$ satisfies the positivity assumption 
(\ref{perpos}).
\end{enumerate}
\end{lm5}

\nd
{\sl Proof.\,\,}  (i) We prove (\ref{A_n exp}) by induction on 
$|\lambda| = \lambda_1 + \cd + \lambda_i$. 
For $|\lambda| = 0$, the sum on the right hand side of 
(\ref{A_n exp}) telescopes to $x_{j;i} - x_{j+i;0} = x_{j;i}$,
as required. 
So we assume that $|\lambda| > 0$.
Let $k \in [1,i]$ be the maximal index such that $\lambda_k >0$,
and let 
$\lambda' = (\lambda_1, \cd, \lambda_{k-1}, \lambda_k -1, 0, \cd, 0)$. 
By (\ref{S-strings}) and (\ref{S-forms A_n}), 
$$
\vp^{(\lambda)}_{j;i}= S_{j+k-1;i-k+\lambda_k} \vp^{(\lambda')}_{j;i}.
$$
By the inductive assumption, (\ref{A_n exp}) holds for 
$\vp^{(\lambda')}_{j;i}$; in particular, the coefficient of 
$x_{j+k-1;i-k+\lambda_k}$ in $\vp^{(\lambda')}_{j;i}$ is equal to $1$.
By (\ref{S_k-periodic}),
\beq
\vp^{(\lambda)}_{j;i} = 
\vp^{(\lambda')}_{j;i} - \beta_{j+k-1;i-k+\lambda_k},
\label{induc A_n}
\eeq
where the forms $\beta_{j;i}$ are given by
\beq
\beta_{j;i}(\vec x) = x_{j;i}-x_{j;i+1}- x_{j+1;i-1}+x_{j+1;i}.
\label{beta A_n}
\eeq
Here note that 
$x_{j;i+1}=0$ if $i=n$ and $x_{j+1;i-1}=0$ if $i=1$.
Expressing the two summands on the right hand side of (\ref{induc A_n})
via (\ref{A_n exp}) and (\ref{beta A_n}) respectively, 
we see that (\ref{A_n exp}) is also valid for $\vp^{(\lambda)}_{j;i}$.
This completes the proof of (i).

(ii) This is just a special case of (\ref{A_n exp}).

(iii) We only need to show that the set of forms
$\vp^{(\lambda)}_{j;i}$
is closed under all the transformations $S_{j';i'}$.
An easy check using (\ref{A_n exp}) and (\ref{beta A_n}) shows 
that the action of  
$S_{j';i'}$ on  $\vp^{(\lambda)}_{j;i}$ can be described as follows. 
For an $i$-admissible partition $\lm=(\lm_1,\cd,\lm_i)$ and 
$k = 1, \cd, i$, we denote by 
$\lm\leftarrow k$ (resp. by $\lm\rightarrow k$) 
the sequence obtained from $\lm$ by replacing $\lm_k$ with $\lm_k + 1$
(resp. with $\lm_k -1$) provided that this sequence is an 
$i$-admissible partition; otherwise, we set $\lm\leftarrow k = \lm$ 
(resp. $\lm\rightarrow k = \lm$). 
Then we have
\beq
S_{j';i'}\vp^{(\lm)}_{j;i}=
\cases{
\vp^{(\lm\leftarrow k)}_{j;i}& if $(j';i') = (j+k-1;i-k+1+\lm_k)$,\cr
\vp^{(\lm\rightarrow k)}_{j;i}& if $(j';i') = (j+k;i-k+\lm_k)$,\cr
\vp^{(\lm)}_{j;i}& otherwise.\cr}
\label{non-triv S}
\eeq

(iv) In view of (i) and (iii), it is enough to observe that 
the only components that can occur in $\vp^{(\lambda)}_{j;i}$
with a negative coefficient are $x_{j+k;i-k+\lambda_k}$ for some $k \geq 1$.
Since $j+k \geq 2$, (\ref{perpos}) follows. 
This completes the proof of the lemma. \qed

\vskip5pt

Now we can complete the proof of Theorem \ref{A_n}.
Using Theorem \ref{main} and parts (iii) and (iv) of Lemma \ref{lm5.1},
we conclude that
${\rm Im} \,(\Psi_{\io})$ is the set of all nonnegative integer families
$(x_{j;i})$ such that $\vp^{(\lambda)}_{j;i} \geq 0$ for all
$(j;i) \in \ZZ_{\geq 1} \times [1,n]$ and all 
$i$-admissible partitions $\lambda$.
If $i = 1$, and $\lambda =(l)$ is a $1$-admissible partition
(i.e., $l \in [1,n]$) then (\ref{A_n exp}) gives
\beq
\vp^{(l)}_{j;1} = x_{j;l+1} - x_{j+1;l}.
\label{j;1}
\eeq
Combining the inequalities $\vp^{(l)}_{j;1} \geq 0$
with the inequalities $x_{j+i;n+1-i} \leq 0$ provided
by Lemma \ref{lm5.1} (ii), we obtain the desired set of inequalities 
in Theorem \ref{A_n}.
It remains to show that all the inequalities 
$\vp^{(\lambda)}_{j;i} \geq 0$ are consequences 
of the ones with $i=1$. 
But this follows at once from (\ref{A_n exp}) and (\ref{j;1}),
which can be written as
$$\vp^{(\lambda)}_{j;i} = 
\sum_{k=1}^i \vp^{(i-k+\lambda_k)}_{j+k-1;1}.$$
Theorem \ref{A_n} is proved. \qed

\section{$A^{(1)}_{n-1}$-case}
\setcounter{equation}{0}
\renewcommand{\theequation}{\thesection.\arabic{equation}}
In this section we shall apply 
Theorem \ref{main} to the case when $\ge$ is an affine Lie algebra
of type $A^{(1)}_{n-1}$ (also sometimes denoted by  
$\widehat {\ssl_{n}}$). 
We will assume that $n \geq 3$ since the case of $A^{(1)}_{1}$
was already treated above. 
We will identify the index set $I$ with $[1,n]$
in the standard way; thus, the Cartan matrix 
$(a_{i,j}= \lan h_i,\al_j\ran )_{1 \leq i,j \leq n}$ is given by 
$a_{i,i}=2$, $a_{i,j}=-1$ for $|i-j|=1$ or $|i-j|=n-1$, and 
$a_{i,j}=0$ otherwise. 
We will find the image of the Kashiwara embedding 
${\rm Im} \,(\Psi_{\io}) \subset \ZZ^{\ify}_{\geq 0}$
for the periodic sequence 
$$
\io = (\cd,\underbrace{n,\cd,2,1}_{},
\cd,\underbrace{n,\cd,2,1}_{},\underbrace{n,\cd,2,1}_{}).
$$

To formulate the answer, we need some terminology and notation.
For any $k \geq 1$, let $\Xi_k = \Xi_{k,\io}$ denote the set of forms
that can be obtained from $x_k$ by a sequence of piecewise-linear 
transformations $S_j$ (cf (\ref{Xi})). 
In dealing with $\Xi_k$, we will use the shorthand
$$j;i[k] := k-1 + (j-1)(n-1) + i.$$
Thus, the correspondence $(j;i) \mapsto j;i[k]$ is a bijection between
$\ZZ_{\geq 1} \times [1,n-1]$ and $\ZZ_{\geq k}$. 
This bijection transforms the usual linear order on $\ZZ_{\geq k}$
into the {\it lexicographic} order on $\ZZ_{\geq 1} \times [1,n-1]$
given by
$$(j';i') < (j;i) \,\,{\rm if}\,\, j' < j \,\,{\rm or} \,\, 
j'=j,\, i' < i.$$
We will consider integer ``matrices"  $C = (c_{j;i})$ indexed by
$\ZZ_{\geq 1} \times [1,n-1]$, and such that $c_{j;i} = 0$
for $j \gg 0$.
With every such $C$ and any $k \geq 1$ we associate a linear form
$\vp_{C[k]}$ on $\ZZ^{\ify}$ given by
\beq
\vp_{C[k]} = \sum_{j;i} c_{j;i} x_{j;i[k]}.
\label{forms-affine A}
\eeq 
For any $(j;i) \in \ZZ_{\geq 1} \times [1,n-1]$, we set
$$s_{j;i} = s_{j;i}(C) = c_{1;i} + c_{2;i} + \cd + c_{j;i}.$$
We will say that a matrix $C$ 
(and each of the corresponding forms $\vp_{C[k]}$)
 is {\it admissible} if it satisfies the following conditions:
\beq
s_{j;i} \geq 0 \,\,{\rm for}\,\, (j;i) \in \ZZ_{\geq 1} \times [1,n-1].
\label{adm1}
\eeq
\beq
s_{j;i} = \delta_{i,1} \,\,{\rm for}\,\, j \gg 0.
\label{adm2}
\eeq
\beq
\sum_{(j';i') \leq (j;i)} s_{j';i'} \leq j
\,\,{\rm for \,\, any}\,\, (j;i), 
\,\,{\rm with \,\, the \,\, equality \,\, for}
\,\, j \gg 0.
\label{adm3}
\eeq
\beq
{\rm If}\,\, s_{j;i} > 0 \,\,{\rm then}\,\, s_{j';i'} > 0
\,\,{\rm for \,\, some}\,\, (j';i') \,\,{\rm with}\,\, 
(j;i) < (j';i') \leq (j+1;i).
\label{adm4}
\eeq

As an example, fix $(j;i)$ and take $C$ such that 
the only non-zero terms $s_{j';i'}$ are $s_{j;i} = j$
and $s_{j';1} = 1$ for all $j' > j$.
The admissibility conditions are obviously satisfied, and
the corresponding admissible forms $\vp_{C[k]}$ are given by
$$\vp_{C[k]} = j x_{j;i[k]} + x_{j+1;1[k]} - j x_{j+1;i[k]}.$$
In particular, if $C_0$ is the matrix corresponding to $(j;i) = (1;1)$ 
then $\vp_{C_0 [k]} = x_k$.

\newtheorem{thm6}{Theorem}[section] 
\begin{thm6}
\label{affine A}
The image ${\rm Im} \,(\Psi_{\io})$
of the Kashiwara embedding is the set of 
$\vec x \in \ZZ^{\ify}$ such that 
$\vp_{C[k]} (\vec x) \geq 0$ for all admissible forms $\vp_{C[k]}$. 
\end{thm6}

\vskip5pt
\nd
{\sl Proof.\,\,}
The following lemma constitutes the main part of the proof.

\newtheorem{lm6}[thm6]{Lemma}
\begin{lm6}
\label{lm1 aff}
For any $k \geq 1$, the set $\Xi_k = \Xi_{k,\io}$ of forms
that can be obtained from $x_k$ by a sequence of piecewise-linear 
transformations $S_j$, consists of all admissible forms $\vp_{C[k]}$.
\end{lm6}

Before proving this lemma, let us show that it implies our theorem.
In view of Theorem \ref{main}, it suffices to show that
every admissible form satisfies the positivity assumption (\ref{perpos}). 
In other words, we need to show that, for every admissible matrix $C$,
all the entries $c_{1;i}$ ($1\leq i\leq n-1$) and $c_{2;1}$ are nonnegative. 
By (\ref{adm1}), we have $c_{1;i} = s_{1;i} \geq 0$;
so it remains to show that $c_{2;1} \geq 0$. 
Again using (\ref{adm1}), if we assume $c_{2;1} < 0$ then we must have
$c_{1;1} = s_{2;1} - c_{2;1} > 0$. 
The proof of (\ref{perpos}) is now completed by the following lemma.

\begin{lm6}
\label{lm2 aff}
The matrix $C_0$ with the entries 
$c_{1;1} = 1$, and $c_{j;i} = 0$ for $(j;i) \neq (1;1)$, is
the only admissible matrix with $c_{1;1} > 0$.
\end{lm6}
{\sl Proof.\,\,}
Combining the condition $c_{1;1} > 0$ with (\ref{adm1}) and
(\ref{adm3}) (for $(j;i) = (1;n-1)$), 
we conclude that $s_{1;1} = c_{1;1} = 1$, 
and $s_{1;i} = 0$ for $i \neq 1$.
Now (\ref{adm4}) implies that $s_{2;1} > 0$.
Combining the latter condition with (\ref{adm1}) and
(\ref{adm3}) (for $(j;i) =(2;n-1)$), we conclude that $s_{2;1} = 1$, 
and $s_{2;i} = 0$ for $i \neq 1$.
Continuing in the same manner, we conclude that 
$s_{j;i} = \delta_{i,1}$ for all $(j;i)$, i.e., $C = C_0$. \qed

\vskip5pt

We now turn to the proof of Lemma \ref{lm1 aff}.
First let us show that, for any $k \geq 1$, the set of admissible 
forms $\vp_{C[k]}$ is closed under all transformations
$S_{j;i[k]}$. 
Note that the forms $\beta_{j;i[k]}$ are given by
\beq
\beta_{j;i[k]} = x_{j;i[k]} - x_{j;i+1[k]} - x_{j+1;i[k]}
+ x_{j+1;i+1[k]},
\label{beta aff}
\eeq
with the convention that $x_{j;n[k]} = x_{j+1;1[k]}$
The definitions readily imply that $S_{j;i[k]} \vp_{C[k]} = \vp_{C'[k]}$,
where the matrix $C'$ is obtained from $C$ as follows.
If $c_{j;i} = 0$ then $C' = C$.
Otherwise, we set $s_{j';i'} = s_{j';i'} (C)$ and 
$s'_{j';i'} = s_{j';i'} (C')$.
Then the only values $s'_{j';i'}$ that are different from $s_{j';i'}$
are given by:

(I) If $c_{j;i} > 0$ then 
$s'_{j;i} = s_{j;i} - c_{j;i} = s_{j-1;i}$ and  
$s'_{j;i+1} = s_{j;i+1} + c_{j;i}$. 

(II) If $c_{j;i} < 0$ then 
$s'_{j-1;i-1} = s_{j-1;i-1} - c_{j;i}$ and  
$s'_{j-1;i} = s_{j-1;i} + c_{j;i} = s_{j;i}$.

\nd
Note that in case (II), if $i = 1$ then (as shown above) we must have
$j > 2$, and we identify  the index $(j-1;i-1)$ with $(j-2;n-1)$.
We need to show that in both cases, the transformation $C \mapsto C'$
preserves admissibility.

In both cases, the conditions (\ref{adm1}) and (\ref{adm2}) 
for $C'$ are obvious.
To prove that $C'$ satisfies (\ref{adm4}), we notice
that, in case (I),   
$$s'_{j;i+1} = s_{j;i+1} + c_{j;i} > 0, \,\,
s_{j;i} = s_{j;i-1} + c_{j;i} > 0;$$
similarly, in case (II),
$$s'_{j-1;i-1} = s_{j-1;i-1} - c_{j;i} > 0, \,\,
s_{j-1;i} = s_{j;i} - c_{j;i} > 0.$$
The fact that $C'$ satisfies (\ref{adm4}), is now a consequence of 
the following lemma.

\begin{lm6}
\label{lm3 aff}
Suppose the nonnegative integer families $s = (s_{j;i})$ and 
$s' = (s'_{j;i})$ satisfy one of the following two sets of properties:
\begin{enumerate}
\item[{\rm(i)}] For some $(j;i)$, we have $s'_{j;i} = s_{j-1;i}$,
$s_{j;i} >0$, $s'_{j;i+1} > 0$, and $s'_{j';i'} = s_{j';i'}$
for $(j';i')$ different from $(j;i)$ and $(j;i+1)$.
\item[{\rm(ii)}] For some $(j;i)$, we have $s'_{j-1;i} = s_{j;i}$,
$s_{j-1;i} >0$, $s'_{j-1;i-1} > 0$, and $s'_{j';i'} = s_{j';i'}$
for $(j';i')$ different from $(j-1;i-1)$ and $(j-1;i)$.
\end{enumerate}
Then if $s$ satisfies (\ref{adm4}), the same is true for $s'$.
\end{lm6}
{\sl Proof.\,\,}
We will prove the statement in case (i), the argument in case (ii) being
very similar.  
The failure of (\ref{adm4}) for $s'$ means that,
for some $(j_0;i_0)$, we have $s'_{j_0;i_0} > 0$ 
and $s'_{j';i'} = 0$ for $(j_0;i_0) < (j';i') \leq (j_0+1;i_0)$;
we will refer to this as the $(j_0;i_0)$-violation.
Since passing from $s$ to $s'$ only affects the entries
$(j;i)$ and $(j;i+1)$, the $(j_0;i_0)$-violation for $s'$ 
can only occur when $(j-1;i) \leq (j_0;i_0)\leq (j;i+1)$. 
Furthermore, since $s'_{j;i+1} > 0$, there is no $(j_0;i_0)$-violation
for $(j-1;i) < (j_0;i_0) < (j;i+1)$.
The $(j-1;i)$-violation for $s'$ is impossible because 
$s'_{j;i} = s_{j-1;i} = s'_{j-1;i}$.
Finally, since $s_{j;i} > 0$, the $(j;i+1)$-violation for $s'$
would imply the $(j;i+1)$-violation or the $(j;i)$-violation for $s$,
hence is also impossible. \qed

\vskip5pt

Continuing the proof of Lemma \ref{lm1 aff},
the fact that $C'$ satisfies the remaining condition (\ref{adm3})
is obvious in case (I). 
In case (II), the failure of (\ref{adm3}) for $C'$ 
can only happen when $C \mapsto C'$ affects the entries
$s_{j-1;0} = s_{j-2;n-1}$ and $s_{j-1;1}$ for some $j > 2$;
this possible violation of (\ref{adm3}) is
$$\sum_{(j';i') < (j-1;1)} s'_{j';i'} \geq j-1.$$
However, since $C$ satisfies (\ref{adm3}), we also have
$$\sum_{(j'; i') < (j;1)} s'_{j';i'} =
\sum_{(j'; i') < (j;1)} s_{j';i'} \leq j-1.$$
Combining the last two inequalities and using (\ref{adm1}),
we conclude that $s'_{j-1;i'} = 0$ for all $i' \in [1,n-1]$.
But this contradicts the fact that $C'$ 
satisfies (\ref{adm4}), which we already established.

To complete the proof of Lemma \ref{lm1 aff}, it is enough to show that
every admissible matrix $C'$ can be obtained from the matrix $C_0$
by a sequence of transformations of type (I) above. 
We introduce the linear order on the set of 
admissible matrices by setting $C \prec C'$ 
if $s_{j;i} > s'_{j;i}$ for the minimal index $(j;i)$ 
such that $s_{j;i} \neq s'_{j;i}$. 
In view of Lemma \ref{lm2 aff}, $C_0$ is minimal with respect 
to this order.
Note also that if $C \mapsto C'$ is a transformation of type (I)
then $C \prec C'$. 

Now let $C'$ be an admissible matrix different from $C_0$.
We will construct an admissible matrix $C$ such that 
$C \mapsto C'$ is a transformation of type (I).
By Lemma \ref{lm2 aff}, we have $s'_{1;1} = c'_{1;1} = 0$. 
Let $(j_0;i) \in \ZZ_{\geq 1} \times [1,n-1]$ be the minimal index 
such that $s'_{j_0;i+1} > 0$ (the existence of $(j_0;i)$
is guaranteed by  (\ref{adm2})).
We claim that, for some $j \geq j_0$, 
\beq
\sum_{(j;i) < (j';i') \leq (j+1;i)} s'_{j';i'} > 1.
\label{exist j}
\eeq
Indeed, assuming that (\ref{exist j}) is false for all $j \geq j_0$,
we would obtain
$$\sum_{(j';i') \leq (j;i)} s'_{j';i'} \leq j - j_0 < j$$
for $j \gg 0$, which contradicts (\ref{adm3}). 
Let $j \geq j_0$ be the minimal index satisfying (\ref{exist j}).
Arguing as in the proof of Lemma \ref{lm2 aff}, we conclude 
that the only non-zero terms $s'_{j';i'}$ for
$(j_0;i) \leq (j';i') \leq (j;i)$ are
$s'_{j';i+1} = 1$ for $j_0 \leq j' < j$. 
Furthermore, we have $s'_{j;i+1} > 0$, and if $s'_{j;i+1} = 1$
then $s'_{j';i'} > 0$ for some $(j';i')$ with 
$(j;i+1) < (j';i') \leq (j+1;i)$. 
Now we define the matrix $C$ by setting 
$s_{j;i} = s_{j;i} (C) = s'_{j;i} + 1$, 
$s_{j;i+1} = s'_{j;i+1} - 1$, and $s_{j';i'} = s'_{j';i'}$
for $(j';i')$ different from $(j;i)$ and from $(j;i+1)$. 
The definitions readily imply that $C$ is admissible and 
$C \prec C'$, and 
that $S_{j;i}:C \mapsto C'$ is a transformation of type (I). 
Iterating this construction, we see that
$C'$ can be obtained from $C_0$
by a sequence of transformations of type (I).
This completes the proof of Lemma \ref{lm1 aff} and then
Theorem \ref{affine A}. \qed

\vskip5pt
\nd
{\bf Remarks.} \,\, (a)	A direct check using (\ref{beta aff}) 
shows that the form $\vp_{C[k]}$ can also be written as
\beq
\vp_{C[k]} = x_k - \sum_{(j;i)\geq (1;1)} d_{j;i[k]} \beta_{j;i[k]},
\label{beta expression}
\eeq
where the coefficients $d_{j;i[k]}$ are given by
$$d_{j;i[k]} = j - \sum_{(j';i') \leq (j;i)} s_{j';i'}.$$
Thus, the meaning of (\ref{adm3}) is that the sum in (\ref{beta expression})
is a (finite) nonnegative linear combination of the  
$\beta_{j;i[k]}$.

\vskip5pt

(b) It would be interesting to find the minimal set 
of inequalities defining ${\rm Im} \,(\Psi_{\io})$, i.e., to eliminate
the redundant linear forms among the $\vp_{C[k]}$.


\begin{thebibliography}{99}


\def\CMP{\sl Commum.Math.Phys.}
\def\IJMP{\sl Int.J.Mod.Phys.}
\def\Duke{\sl Duke Math.J.}

\bibitem{BFZ} Berenstein A, Fomin S and Zelevinsky A, 
Parametrizations 
of canonical bases and totally positive matrices,
{\sl Adv. \ in \ Math.},  {\bf 122} (1996), 49--149.

\bibitem{BZ1} Berenstein A and Zelevinsky A, 
String bases for quantum groups of type $A_r$, 
{\sl Advances in Soviet Math.}, {\bf 16}, Part 1 (1993),
51--89.

\bibitem{BZ2} Berenstein A and Zelevinsky A, 
Canonical bases for the quantum group
of type $A_r$ and piecewise-linear combinatorics, 
{\sl Duke Math J.}, {\bf 82} (1996), 473-502.

\bibitem{CP} Chari V and Pressley A, A guide to Quantum Groups,
	    Cambridge Univ.Press (1994).

\bibitem{JMMO} Jimbo M, Misra K.C, Miwa T and Okado M,
Combinatorics of representatins of $U_q(\what{\ssl(n)})$ at $q=0$,
{\sl Comm. Math. Phys.}, {\bf 136} (1991), 543--566.
\bibitem{K0} Kashiwara M, Crystallizing the $q$-analogue of 
	   universal enveloping algebras, Comm. Math. Phys., {\bf 133},
	   (1990), 249--260.
\bibitem{K1} Kashiwara M,
 On crystal bases of the $q$-analogue of universal enveloping algebras,
	{\it Duke Math. J.},{\bf 63} (1991) 465--516.


\bibitem{K3} Kashiwara M,
      Crystal base and Littelmann's refined Demazure character formula,
	     {\sl Duke Math. J.}, {\bf 71}(1993), 839--858.

\bibitem{K4} Kashiwara M,
      Crystal base of modified quantized enveloping algebra,
	     {\sl Duke Math. J.}, {\bf 73} (1994), 383--413  .



\bibitem{KMN1}
       Kang S-J, Kashiwara M, Misra K, Miwa T, Nakashima T and Nakayashiki A,
  Affine crystals and vertex models, 
  {\IJMP},{A7 }Suppl. 1A (1992) 449--484.
\bibitem{KMN2}
       Kang S-J, Kashiwara M, Misra K, Miwa T, Nakashima T and Nakayashiki A,
  Perfect crystals of quantum affine Lie algberas, 
    {\sl Duke Math. J.}, {\bf 68} (1992), 499-607.

\bibitem{KN} Kashiwara M and Nakashima T, 
	     Crystal graph for representations 
	     of the $q$-analogue of classical Lie algebras, 
	    {\sl J. Algebra}, {\bf 165}, (1994), 295--345.
\bibitem{Kac} Kac V.G , Infinite dimensional Lie algebras, 
	      3rd edition, Cambridge Univ. Press (1990).
\bibitem{Kas} Kassel C, Quantum Groups,
	     GTM 155, Springer-Verlag, (1995).
\bibitem{L1} Littelmann P, A Littlewood-Richardson type rule for 
	  symmetrizable Kac-Moody algebras, {\sl Invent. Math.},
	  {\bf 116} (1994), 329--346.
\bibitem{L2} Littelmann P, Path and root operators in representation theory,
	    {\sl Ann. of Math.}, (to appear).

\bibitem{Lu} Lusztig G, {\it Introduction to quantum groups},
Birkh\"auser, Boston, 1993.
\bibitem{N} Nakashima T, Crystal Base and a Generalization of 
the Littlewood-Richardson Rule for the Classical Lie Algebras, 
{\sl Commun. Math. Phys.}, {\bf 154}, (1993), 215-243.
\end{thebibliography}
\end{document}